\documentclass[12pt,reqno]{amsart}
\usepackage{amssymb,amsthm,amsmath,amstext,amsxtra}
\usepackage{mathrsfs,bm}
\usepackage{fullpage}
\usepackage[all]{xy}
\usepackage{booktabs}
\usepackage{mathtools}
\usepackage{hyperref}
\hypersetup{colorlinks=true,urlcolor=blue,citecolor=blue,linkcolor=blue}
\usepackage{enumerate}
\usepackage{ stmaryrd }
\usepackage{enumitem}
\usepackage{comment}
\usepackage{multirow}
\usepackage{colonequals}
\usepackage{tikz}
\usepackage{marginnote}
\usepackage[export]{adjustbox}

\usepackage{tikz}
\usepackage{tikz-cd}
\usepackage{xcolor}

\numberwithin{equation}{section}

\begin{document}

\title{Ap\'ery's irrationality proof, Beukers's modular forms and mirror symmetry}

\author{Wenzhe Yang}
\address{SITP, Physics Department, Stanford University, CA, 94305}
\email{yangwz@stanford.edu}

\begin{abstract}
In this paper, we will apply the ideas from the mirror symmetry of Calabi-Yau threefolds to study the modular forms and one-parameter family of K3 surfaces found by Beukers and Peters, which provide enlightenment to the two mysterious sequences constructed by Ap\'ery in his proof of the irrationality of $\zeta(3)$. We will construct a fourth order differential operator and a prepotential from the canonical solutions of this differential operator. The third derivative of this prepotential with respect to the mirror map defines a Yukawa coupling that is a weight-4 modular form. The instanton expansion of this Yukawa coupling yields integral instanton numbers, which are also periodic with period 6.
\end{abstract}

\maketitle
\setcounter{tocdepth}{1}
\vspace{-13pt}
\tableofcontents
\vspace{-13pt}

\section{Introduction}

In 1978, Ap\'ery proved the irrationality of $\zeta(3)$ in a rather miraculous way \cite{Apery}. He constructed two mysterious sequences $\{A_n\}_{n \geq 0}$ and $\{B_n\}_{n \geq 0}$, both of which are solutions for the same recursion equation with different initial conditions. The quotient $B_n/A_n$ converges to $\zeta(3)$ fast enough so that Dirichlet's irrationality criterion can be used to show $\zeta(3)$ is irrational. The two sequences $A_n$ and $B_n$ do have a very ad hoc appearance. Since then, there have been many papers which try to make sense of them. Most notably, Beukers and his collaborators have shown that the two sequences are closely related to the theories of modular forms and the Picard-Fuchs equation of a one-parameter family of K3 surfaces. In this paper, we will borrow ideas from the mirror symmetry of Calabi-Yau threefolds to provide further enlightenment to the works of Ap\'ery, Beukers and many others. In the process, we will obtain a more meaningful prepotential for the one-parameter family of K3 surfaces, whose third derivative with respect to the mirror map gives us a nontrivial Yukawa coupling that is a weight-4 modular form. Similarly, we can also define the instanton numbers for this Yukawa coupling, which are integral and periodic. Before we discuss our results, let us first briefly review the works of Ap\'ery, Beukers and Peters.

\subsection{Beukers's modular forms and K3 surfaces} \label{sec:Beukersmodularforms}

The two mysterious sequences $\{A_n\}_{n \geq 0}$ and $\{B_n\}_{n \geq 0}$ obtained by Ap\'ery are the solutions for the recursion equation
\begin{equation} \label{eq:aperyrecursion}
(n+1)^3x_{n+1}-(34n^3+51n^2+27n+5)x_n+n^3x_{n-1}=0
\end{equation}
with initial conditions
\begin{equation}
(A_0,A_1)=(1,5)~\text{and}~(B_0,B_1)=(0,6).
\end{equation}
Ap\'ery proved that $A_n \in \mathbb{Z}$ for all $n \geq 0$, which is a very remarkable property. Since if we compute $A_n$ recursively via \eqref{eq:aperyrecursion}, we need to divide by $n^3$ at each step, therefore we should a priori expect that $A_n$ has denominator $(n!)^3$. Hence the integrality of $A_n$ is very surprising. Equally surprising is the property that there is a very good denominator bound for $B_n$
\begin{equation}
d_n^3 B_n \in \mathbb{Z},~\forall n \geq 0, ~\text{where}~ d_n=\text{l.c.m.} \{1,2,\cdots,n \}.
\end{equation}
These two properties are crucial in Ap\'ery's proof of the irrationality of $\zeta(3)$ \cite{Apery, Zagier}.

In the paper \cite{Beukers}, Beukers found two modular forms for the modular group $\Gamma_1(6)$ with their $q$-expansions of the form
\begin{equation} \label{eq:Beukerstwomodularforms}
\begin{aligned}
T(q)&=q \prod_{n=1}^\infty \frac{(1-q^n)^{12}(1-q^{6n})^{12}}{(1-q^{2n})^{12}(1-q^{3n})^{12}}=q-12 q^2+66 q^3-220 q^4+\cdots, \\
F(q)&=\prod_{n=1}^\infty \frac{(1-q^{2n})^{7}(1-q^{3n})^{7}}{(1-q^{n})^{5}(1-q^{6n})^{5}}=1+5q+13q^2+23q^3+\cdots.
\end{aligned}
\end{equation}
Here the variable $q$ is by definition
\begin{equation}
q:=\exp 2 \pi i \,\tau,~\tau \in \mathbb{H},
\end{equation}
where $\mathbb{H}$ is the upper-half complex plane. The modular form $F(q)$ has an expansion with respect to $T(q)$ of the form
\begin{equation} \label{eq:modularformsintegrality}
F(q)=\sum_{n=0}^{\infty} A_n T(q)^n,
\end{equation}
which explains the integrality of $A_n$ since both modular forms have integral coefficients. 

In 1984, Beukers and Peters discovered a very interesting relation between the two modular forms and a one-parameter family of K3 surfaces \cite{BeukersPeters}. They found that for a generic $\varphi$, the surface 
\begin{equation}
V_\varphi:~1-\varphi \frac{x(1-x)y(1-y)z(1-z)}{1-z+xyz} =0
\end{equation}
is birationally equivalent to a K3 surface with Picard number 19. Denote this one-parameter family of K3 surfaces by 
\begin{equation} \label{eq:introK3family}
\pi:\mathscr{V} \rightarrow \mathbb{P}^1.
\end{equation}
The nowhere-vanishing holomorphic twoform $\omega_\varphi$ satisfies a Picard-Fuchs equation
\begin{equation}
\mathcal{L}\,\omega_\varphi =0,
\end{equation}
where the Picard-Fuchs operator $\mathcal{L}$ is given by
\begin{equation} \label{eq:introPFoperatorK3}
\mathcal{L}=\vartheta^3-\varphi(34 \vartheta^3+51\vartheta^2+27 \vartheta+5)+\varphi^2(\vartheta+1)^3,~\vartheta=\varphi \frac{d}{d \varphi}.
\end{equation}
A fundamental solution of $\mathcal{L}$ is the power series
\begin{equation} \label{eq:k3aphiseries}
\varpi_0(\varphi):=\sum_{n=0}^\infty A_n \varphi^n=1+5\varphi+73 \varphi^2+1445 \varphi^3+\cdots,
\end{equation}
The exist two other canonical solution for $\mathcal{L}$ that are of the form
\begin{equation} \label{eq:introOthertwosolutions}
\begin{aligned}
\varpi_1(\varphi)&=\frac{1}{2 \pi i}\left( \varpi_0(\varphi) \log \varphi +h_1(\varphi) \right)\\
\varpi_2(\varphi)&=\frac{1}{(2 \pi i)^2}\left( \varpi_0(\varphi) \log^2 \varphi +2 h_1(\varphi) \log \varphi+h_2(\varphi) \right),\\
\end{aligned}
\end{equation}
where $h_1$ and $h_2$ are power series that satisfy $h_1(0)=h_2(0)=0$. Now, let $\tau$ be the quotient 
\begin{equation} \label{eq:introMirrorMap}
\tau=\varpi_1(\varphi)/\varpi_0(\varphi),
\end{equation}
which is also called the mirror map for the family \eqref{eq:introK3family}. The mirror map can be inverted order by order, which gives us the $q$-expansion of $\varphi$. Beukers has shown that $\varphi(q)$ is just the modular form $T(q)$ \cite{Beukers, Zagier}
\begin{equation} \label{eq:introductionmirrormap}
\varphi(q)=T(q).
\end{equation}
Furthermore, the Picard-Fuchs operator $\mathcal{L}$ is the symmetric square of a second order differential operator. The three solutions $\{\varpi_i \}_{i=0}^2$ satisfy the equation
\begin{equation}
\varpi_0 \varpi_2=\varpi_1^2.
\end{equation}

\subsection{The prepotential and Yukawa coupling}

In this paper, we will construct a power series $\Pi_0(\phi)$ from the sequence $ \{ B_n \}_{n \geq 0}$ by
\begin{equation}
\Pi_0(\phi)=\sum_{n=0}^{\infty}B_n \phi^{n+1}=6\phi^2+\frac{351}{4}\phi^3+\frac{62531}{36} \phi^4+\frac{11424695}{288} \phi^5 + \cdots.
\end{equation}
It is a solution of a fourth order differential equation $\mathcal{D} \Pi_0(\phi)=0$. There is a surprising relation between the two operators $\mathcal{D}$ and $\mathcal{L}$. Namely, if we define a transformation between the variables $\phi$ and $\varphi$ by
\begin{equation}
\varphi =1/ \phi,
\end{equation}
then a straightforward computation shows that
\begin{equation} \label{eq:introDoperator}
\mathcal{D}=\varphi^{-2} \vartheta \cdot \mathcal{L}.
\end{equation}
Hence the three independent solutions $\{\varpi_i \}_{i=0}^2$ of $\mathcal{L}$ are also solutions of the fourth order differential operator $\mathcal{D}$.

What is more, in a small neighborhood of $\varphi=0$, there exists a fourth canonical solution for $\mathcal{D}$ of the form
\begin{equation} \label{eq:introvarpi3}
\varpi_3(\varphi)=\frac{1}{(2 \pi i)^3}\left( \varpi_0(\varphi) \log^3 \varphi+3 h_1(\varphi) \log^2 \varphi+3 h_2(\varphi) \log \varphi+h_3(\varphi) \right),
\end{equation}
The four canonical solutions $\{\varpi_i(\varphi) \}_{i=0}^3$ form a basis for the solution space of $\mathcal{D}$. The first question that comes to our mind is whether the fourth canonical solution $\varpi_3$ can tell us an interesting new? The answer is yes. Inspired by the mirror symmetry of Calabi-Yau threefolds,  we will define a prepotential by
\begin{equation}
\mathscr{F}_{\mathcal{D}}=\varpi_3/\varpi_0.
\end{equation}
The Yukawa coupling $\mathscr{Y}_{\mathcal{D}}$ is by definition the third derivative of this prepotential  with respect to the mirror map $\tau$ \eqref{eq:introMirrorMap}
\begin{equation}
\mathscr{Y}_\mathcal{D}=\frac{d^3\mathscr{F}_\mathcal{D}}{d\tau^3}.
\end{equation}
It is very surprising that the $q$-expansion of the Yukawa coupling $\mathscr{Y}_{\mathcal{D}}$ is found to be the product of two weight-2 modular form
\begin{equation}
\mathscr{Y}_{\mathcal{D}}=6 F(q)H(q).
\end{equation}
Here $F(q)$ is the weight-2 modular form in the formula \eqref{eq:Beukerstwomodularforms} and $H(q)$ is the following weight-2 modular form
\begin{equation}
H(q)=2 \Theta_{\text{hex}}^2(q^2)-\Theta_{\text{hex}}^2(q),
\end{equation}
where $\Theta_{\text{hex}}(q)$ is the theta series of the planar hexagonal lattice $A_2$ \cite{Conway, OEIStheta}
\begin{equation}
\Theta_{\text{hex}}(q)=\theta_3(0,\tau)\theta_3(0,3\tau)+\theta_2(0,\tau)\theta_2(0,3\tau).
\end{equation}
Furthermore, we can also define the instanton numbers for the Yukawa coupling $\mathscr{Y}_{\mathcal{D}}$ by the equation
\begin{equation}
\mathscr{Y}_{\mathcal{D}}=6+\sum_{k=1}^\infty k^3 N_k \frac{q^k}{1-q^k}.
\end{equation}
We have computed the instanton numbers $N_k$ for small $k$, and the first six instanton numbers are given by
\begin{equation}
N_1 = -42, \, N_2=-39, \, N_3=-44, \, N_4=-39, \, N_5 =-42, \,N_6=-34,
\end{equation}
which are indeed integers. Very surprisingly, we have found that the instanton number $N_k$ is periodic
\begin{equation}
N_{k+6} = N_k.
\end{equation}
We have also shown that our previous analysis can be immediately generalized to the Dwork family of K3 surfaces, the mirror family for the quartic K3 surfaces, which yields similar results.

The outline of this paper is as follows. In Section \ref{sec:mirrorsymmetryCYthreefolds}, we will very briefly review several results of the mirror symmetry of Calabi-Yau threefolds that will be needed in this paper. In Section \ref{sec:YukawaK3}, we will study the one-parameter family \eqref{eq:introK3family} of K3 surfaces . In Section \ref{sec:varpi3}, we will study the canonical solutions of the forth order differential operator $\mathcal{D}$ \eqref{eq:introDoperator}, in particular, the fourth canonical solution $\varpi_3$ \eqref{eq:introvarpi3}. In Section \ref{sec:NewYukawa}, we will construct a prepotential and a Yukawa coupling from the fourth canonical solution $\varpi_3$ \eqref{eq:introvarpi3} of $\mathcal{D}$, and show they do have very interesting properties. In Section \ref{sec:Dworkfamily}, we will generalize our results in Section \ref{sec:NewYukawa} to the Dwork family of K3 surfaces. In Section \ref{sec:furtherprospects}, we will discuss several important open questions.

\section{The mirror symmetry of Calabi-Yau threefolds} \label{sec:mirrorsymmetryCYthreefolds}

In this section, we will briefly review the concept of the mirror symmetry of Calabi-Yau threefolds, and we will focus on the one-parameter mirror pairs \cite{PhilipXenia, KimYang}. We will only introduce the results that will be needed later in this paper. The readers who are not familiar with mirror symmetry could consult the books \cite{CoxKatz, MarkGross}.

Given a mirror pair $(X^\vee,X)$ of Calabi-Yau threefolds, here one-parameter means that the Hodge numbers of $X^\vee$ and $X$ satisfy
\begin{equation}
h^{1,1}(X^\vee)=h^{2,1}(X)=1.
\end{equation}
In this case, the complexified K\"ahler moduli space $\mathscr{M}_K(X^\vee)$ of $X^\vee$ has a very simple description \cite{MarkGross,KimYang}
\begin{equation}
\mathscr{M}_K(X^\vee)=(\mathbb{R}+i\, \mathbb{R}_{> 0})/\mathbb{Z}=\mathbb{H}/\mathbb{Z},
\end{equation}
where $\mathbb{H}$ is the upper-half plane of $\mathbb{C}$. Suppose $e$ forms a basis for $H^2(X^\vee,\mathbb{Z})$ (modulo torsion) and it is also in the K\"ahler cone of $X^\vee$, then every point of $\mathscr{M}_K(X^\vee)$ is of the form $e\,t, \,t \in \mathbb{H}$ \cite{MarkGross}. Moreover, the point $e\,t$ is equivalent to $e\,(t+1)$ under the quotient by $\mathbb{Z}$. The variable $t$ is called the flat coordinate of $\mathscr{M}_K(X^\vee)$ \cite{PhilipXenia, CoxKatz, KimYang}.

The prepotential $\mathcal{F}$ for $X^\vee$ admits an expansion near $t=i\, \infty$ of the form
\begin{equation} \label{eq:Prepotential}
\mathcal{F}=-\frac{1}{6}\, Y_{111}\, t^3 -\frac{1}{2}\, Y_{011}\,t^2-\frac{1}{2}\,Y_{001}\, t-\frac{1}{6}\,Y_{000}+\mathcal{F}^{\text{np}},
\end{equation}
where $\mathcal{F}^{\text{np}}$ is the non-perturbative instanton correction that have a $q$-expansion 
\begin{equation}
\mathcal{F}^{\text{np}}=\frac{1}{(2 \pi i)^3}\sum_{n=1}^{\infty} a_n q^n~\text{with}~q=\exp 2 \pi i \,t.
\end{equation}
The coefficient $Y_{111}$ in the prepotential \eqref{eq:Prepotential} is the topological intersection number of $e$ given by 
\begin{equation} \label{eq:y111}
Y_{111}=\int _{X^\vee} e\wedge e \wedge e,
\end{equation}
which is a positive integer \cite{PhilipXenia,CoxKatz, MarkGross}. The number $Y_{000}$ is equal to \cite{PhilipXenia}
\begin{equation} \label{eq:PhysicistsY000}
Y_{000}=-3\, \chi (X^\vee)\, \frac{\zeta(3)}{(2 \pi i)^3},
\end{equation}
where $\chi(X^\vee)$ is the Euler characteristic of $X^\vee$. However, the computations of the other two numbers  $Y_{011}$ and $Y_{001}$ are more tricky \cite{PhilipXenia}. The Yukawa coupling for $X^\vee$ is given by the third derivative of the prepotential $\mathcal{F}$ with respect to the flat coordinate $t$
\begin{equation} \label{eq:YukawaDefnKahler}
\mathcal{Y}=-\frac{d^3\mathcal{F}}{dt^3}.
\end{equation}
The instanton number $N_k$ for $X^\vee$ is defined by the equation
\begin{equation} \label{eq:YukawaKahler}
\mathcal{Y}=Y_{111}+\sum_{k=1}^\infty k^3 N_k \frac{q^k}{1-q^k},
\end{equation}
which is the Gromov-Witten invariants of $X^\vee$. In general, the Yukawa coupling $\mathcal{Y}$ is very difficult to be computed directly, and that is where mirror symmetry comes to the rescue.

The idea of mirror symmetry is that the Yukawa coupling $\mathcal{Y}$ for $X^\vee$ can be computed by the periods of its mirror family, i.e. a deformation of the mirror threefold $X$. More explicitly, suppose the mirror threefold $X$ admits a deformation over $\mathbb{P}^1$
\begin{equation} \label{eq:mirrorfamily}
\pi:\mathscr{X} \rightarrow \mathbb{P}^1,
\end{equation}
where by abuse of notations the coordinate of $\mathbb{P}^1$ is also denoted by $\varphi$. We will further assume that for each smooth fiber $\mathscr{X}_\varphi$, there exists a nowhere-vanishing holomorphic threeform $\Omega_\varphi$ on $\mathscr{X}_\varphi$ \cite{CoxKatz, MarkGross, KimYang}. From the Griffiths transversality, $\Omega_\varphi$ satisfies a fourth order Picard-Fuchs equation. Suppose further the point $\varphi=0$ is the large complex structure limit of the mirror family \eqref{eq:mirrorfamily}, which means that there exist four homological cycles $\{ \beta_i \}_{i=0}^3$ of $H_3(X,\mathbb{C})$ such that the periods of $\Omega_\varphi$, i.e.
\begin{equation}
\pi_i=\int_{\beta_i} \Omega_\varphi,~i=0,1,2,3,
\end{equation}
are of the form
\begin{equation} \label{eq:PeriodsCan}
\begin{aligned}
\pi_0 &= f_0,  \\
\pi_1 &=\frac{1}{2\pi i}\left(f_0 \log \varphi+f_1\right), \\
\pi_2 &=\frac{1}{(2\pi i)^2}\left( f_0 \log^2 \varphi +2\, f_1 \log \varphi + f_2\right), \\
\pi_3 &=\frac{1}{(2 \pi i)^3} \left( f_0 \log^3 \varphi +3 \, f_1 \log^2 \varphi +3\, f_2 \log \varphi +f_3 \right),
\end{aligned} 
\end{equation}
where $\{f_j\}_{j=0}^3$ are power series in $\varphi$ that satisfy
\begin{equation} \label{eq:boundarycondition}
f_0(0)=1,~f_1(0)=f_2(0)=f_3(0)=0.
\end{equation}
Namely, the monodromy matrix of the column vector $(\pi_0, \pi_1, \pi_2,\pi_3)^\top$ is maximally unipotent. The four canonical periods $\pi_i$ can also be computed by solving the Picard-Fuchs equation satisfied by the holomorphic threeform $\Omega_\varphi$, which is usually enormously simpler \cite{PhilipXenia, MarkGross}.

We will further assume that there exists a nonzero constant $\ell$ such that the  two homological cycles $\ell \beta_0$ and $\ell \beta_1$ are integral, i.e. they lie in $H_3(X, \mathbb{Z})$. In fact, this assumption is satisfied by all the known examples of the one-parameter mirror pairs of Calabi-Yau threefolds. At the heart of mirror symmetry is the construction of the mirror map
\begin{equation} \label{eq:mirrormap}
t=\frac{\pi_1}{\pi_0}=\frac{1}{2 \pi i}\,\left( \log \varphi+\frac{f_1}{f_0} \right),
\end{equation}
which identifies the flat coordinate $t$ on the K\"ahler side with the quotient $\pi_1/\pi_0$ on the complex side. The mirror map \eqref{eq:mirrormap} can be inverted order by order, which yields a $q$-expansion of $\varphi$. If we assume the mirror symmetry conjecture, then the prepotential $\mathcal{F}$ \eqref{eq:Prepotential} can be computed by the canonical periods $\pi_i$ \eqref{eq:PeriodsCan} \cite{PhilipXenia, KimYang}. More explicitly, we have
\begin{equation}\label{eq:prepotentialMirror}
\mathcal{F}=-\frac{1}{12} Y_{111} \left(- \frac{\pi_3}{\pi_0} + \frac{\pi_1 \pi_2}{\pi_0^2}\right)-\frac{1}{2}Y_{011}\frac{\pi_1^2}{\pi_0^2}-\frac{1}{2}Y_{001}\frac{\pi_1}{\pi_0}-\frac{1}{6} Y_{000}.
\end{equation}
On the other hand, the Yukawa coupling $\mathcal{Y}$ \eqref{eq:YukawaDefnKahler} can also be expressed as
\begin{equation} \label{eq:Yukawacomplex}
\mathcal{Y}=\frac{1}{\pi_0^2} \left( \frac{1}{2 \pi i}\frac{d\varphi}{dt} \right)^3 \int_X \Omega_\varphi \wedge \frac{d^3}{d\varphi^3}\Omega_\varphi.
\end{equation}
The instanton expansion of $\mathcal{Y}$ \eqref{eq:YukawaKahler} can be obtained by just plugging in the $q$-expansion of $\varphi$. Notice that the overall scale of formula \eqref{eq:Yukawacomplex} should be chosen such that $\mathcal{Y} \rightarrow Y_{111}$ as $t \rightarrow i \infty$. The readers are referred to the books \cite{CoxKatz, MarkGross} for more details about these computations for the quintic Calabi-Yau threefolds.

\section{The Yukawa coupling for the K3 surfaces} \label{sec:YukawaK3}

In this section, we will show that the Yukawa coupling for the one-parameter family \eqref{eq:introK3family} of K3 surfaces is a trivial, i.e. it is a constant. Therefore, unlike the Calabi-Yau threefolds, the Yukawa coupling does not contain interesting information about the K3 surfaces. In fact, this is a general property for one-parameter families of K3 surfaces, as the Gromov-Witten invariants of K3 surfaces are trivial  \cite{Hartmann, Nagura}.

\subsection{The Frobenius method}

In order for this paper to be as self-contained as possible, let us first show how to find the solutions for the Picard-Fuchs operator $\mathcal{L}$ \eqref{eq:introPFoperatorK3} by using the Frobenius method \cite{YangPeriods}. First, try a solution of the form
\begin{equation}
\mathcal{L}\,\left(\varphi^{\epsilon} \, \sum_{k=0}^{\infty} a_k(\epsilon)\, \varphi^k \right)=0,
\end{equation} 
where $\epsilon$ is a formal variable and $a_k(\epsilon)$ is an unknown number. In order for this equation to be satisfied, to the lowest order we must have
\begin{equation} \label{eq:indicial}
\epsilon^3=0,
\end{equation}
which is called the indicial equation for $\mathcal{L}$. We already know the power series $\varpi_0$ \eqref{eq:k3aphiseries} is a solution of $\mathcal{L}$, hence the indicial equation \eqref{eq:indicial} tells us that the other two independent solutions of $\mathcal{L}$, it they exist, must be of the form in the formula \eqref{eq:introOthertwosolutions}. Now, plug the expression 
\begin{equation}
\varpi_1(\varphi)=\frac{1}{2 \pi i} \left( \varpi_0(\varphi) \log \varphi+\sum_{n=1}^\infty c_{1,n} \varphi^n \right),
\end{equation}
into the Picard-Fuchs equation $\mathcal{L}\varpi_1=0$, and we find that the coefficient $c_{1,n}$ is determined by the recursion equation
\begin{equation} \label{eq:CrecursionequationK3}
\begin{aligned}
&(n+1)^3c_{1,n+1}-\left(34 n^3+51 n^2+27 n+5\right)c_{1,n}+n^3c_{1,n-1}\\
&+3 (n+1)^2A_{n+1}-3 \left(34 n^2+34 n+9\right) A_n+3n^2A_{n-1}=0,
\end{aligned}
\end{equation}
with initial condition \label{eq:C1recursionequationK3initial}
\begin{equation}
c_{1,1}=12,~c_{1,2}=210.
\end{equation}
For simplicity, let us define $h_1(\varphi)$ by
\begin{equation} \label{eq:h1defn}
h_1(\varphi)=\sum_{n=1}^\infty c_{1,n} \varphi^n =12 \varphi +210 \varphi^2+4438 \varphi^3+ \cdots,\\
\end{equation}

Similarly, plug the expression
\begin{equation}
\varpi_2(\varphi)=\frac{1}{(2 \pi i)^2}\left( \varpi_0(\varphi) \log^2 \varphi +2 h_1(\varphi) \log \varphi+\sum_{n=1}^\infty c_{2,n} \varphi^n\right).\\
\end{equation}
into the Picard-Fuchs equation $\mathcal{L}\varpi_2=0$, we obtain the recursion equation for $c_{2,n}$
\begin{equation} \label{eq:C2recursionequationK3}
\begin{aligned}
&(n+1)^3c_{2,n+1}-\left(34 n^3+51 n^2+27 n+5\right)c_{2,n}+n^3c_{2,n-1}\\
&+6 (n+1)^2c_{1,n+1}-6 \left(34 n^2+34 n+9\right) c_{1,n}+6n^2c_{1,n-1} \\
&+6(n+1)A_{n+1}-102(2n+1)A_n+6nA_{n-1}=0
\end{aligned}
\end{equation}
with initial condition
\begin{equation}
c_{2,1}=0,~c_{2,2}=144.
\end{equation}
For simplicity, let $h_2(\varphi)$ be
\begin{equation} \label{eq:h2defn}
h_2(\varphi)=\sum_{n=1}^\infty c_{2,n} \varphi^n=144 \varphi^2 +4320 \varphi^3 +118500 \varphi^3 +\cdots.
\end{equation}
The Picard-Fuchs operator $\mathcal{L}$ is the symmetric square of the following second order differential operator \cite{BeukersPeters, Zagier}
\begin{equation} \label{eq:Lroot}
(1-34\varphi+\varphi^2) \vartheta^2-\varphi (17-\varphi) \vartheta -\frac{1}{4} \varphi(10 - \varphi),
\end{equation}
from which we deduce that
\begin{equation} \label{eq:varpi2defn}
\varpi_0 \varpi_2=\varpi^2_1.
\end{equation}

\subsection{The Yukawa coupling}

Following the mirror symmetry of Calabi-Yau threefolds, let us define the Yukawa coupling $\mathcal{Y}_{bp}$ for the family \eqref{eq:introK3family} by
\begin{equation}
\mathcal{Y}_{bp}= \int_{\mathscr{V}_\varphi} \omega_\varphi \wedge \frac{d^2\omega_\varphi}{d\varphi^2},
\end{equation}
From Griffiths transversality, we have
\begin{equation} \label{eq:GriffithsVanishingparings}
\int_{\mathscr{V}_\varphi} \omega_\varphi \wedge \omega_\varphi=\int_{\mathscr{V}_\varphi} \omega_\varphi \wedge \frac{d\omega_\varphi}{d\varphi}=0,
\end{equation}
from which we obtain
\begin{equation} \label{eq:omegaidentities}
\int_{\mathscr{V}_\varphi} \frac{d\omega_\varphi}{d\varphi} \wedge \frac{d\omega_\varphi}{d\varphi}=-\mathcal{Y}_{bp},~ \int_{\mathscr{V}_\varphi} \frac{d\omega_\varphi}{d\varphi} \wedge \frac{d^2\omega_\varphi}{d\varphi^2}=-\frac{1}{2} \frac{d \mathcal{Y}_{bp}}{d\varphi},~\int_{\mathscr{V}_\varphi} \omega_\varphi \wedge \frac{d^3\omega_\varphi}{d\varphi^3}=\frac{3}{2} \frac{d \mathcal{Y}_{bp}}{d\varphi}.
\end{equation}
The Picard-Fuchs operator $\mathcal{L}$ \eqref{eq:introPFoperatorK3} can also be expressed as
\begin{equation} \label{eq:picardfuchsk3dvarphi}
\begin{aligned}
\mathcal{L}=\,& \varphi ^3 \left(1-34 \varphi +\varphi ^2\right) \frac{d^3}{d\varphi^3}+3 \varphi ^2 \left(1 -51 \varphi +2 \varphi ^2 \right) \frac{d^2}{d\varphi^2} \\
&+\varphi  \left(1 -112 \varphi +7 \varphi ^2\right)\frac{d}{d\varphi}-\varphi (5 -\varphi ) .
\end{aligned}
\end{equation}
From formulas \eqref{eq:picardfuchsk3dvarphi} and \eqref{eq:omegaidentities}, we obtain
\begin{equation}
\frac{3}{2}\varphi ^3 \left(1-34 \varphi +\varphi ^2\right) \frac{d \mathcal{Y}_{bp}}{d\varphi}+3 \varphi ^2 \left(1 -51 \varphi +2 \varphi ^2\right) \mathcal{Y}_{bp}=0,
\end{equation}
whose solution is
\begin{equation} \label{eq:yukawacouplingrational}
\mathcal{Y}_{bp}=\frac{C_{bp}}{\varphi ^2 (1-34 \varphi +\varphi ^2)},
\end{equation}
where $C_{bp}$ is a constant.

On the other hand, the canonical periods $\varpi_i, i=0,1,2$ are also given by the integrals of the holomorphic twoform $\omega_\varphi$ over certain homological cycles in $H_2(\mathscr{V}_\varphi,\mathbb{C})$. There exist cohomological cycles $\{ \alpha_i\}_{i=0}^2 \subset H_2(\mathscr{V}_\varphi,\mathbb{C})$ such that
\begin{equation}
\omega_\varphi=\varpi_0 \alpha_0+\varpi_1 \alpha_1-\varpi_2 \alpha_2.
\end{equation}
While from formula \ref{eq:GriffithsVanishingparings}, we can compute the intersection pairings for $(\alpha_0, \alpha_1, \alpha_2)$, which are given by
\begin{equation} \label{eq:pairingrelations}
\alpha_0 \cdot \alpha_2=\alpha_2 \cdot \alpha_0=C_0,~ \alpha_1 \cdot \alpha_1=2 C_0,~
 \alpha_i \cdot \alpha_j=0 ~\text{otherwise},
\end{equation}
where $C_0$ is a nonzero constant. Hence from formula \ref{eq:omegaidentities}, $\mathcal{Y}_{bp}$ can also be expressed as 
\begin{equation} \label{eq:yukawaperiodsrepn}
\mathcal{Y}_{bp}=-2 C_0 \left( \varpi'^2_1 -\varpi'_0 \varpi'_2 \right)
\end{equation}
Plug in the power series expansion of $\varpi_i$, we find the value of the constant $C_{bp}$
\begin{equation}
C_{bp}=-\frac{2C_0}{(2 \pi i)^2}.
\end{equation}

Following the mirror symmetry of Calabi-Yau threefolds, the mirror map for the one-parameter family \eqref{eq:introK3family} of K3 surfaces is defined by
\begin{equation} \label{eq:k3defnofmirrormap}
\tau=\frac{\varpi_1}{\varpi_0}.
\end{equation}
By abuse of notations, let $q$ also be
\begin{equation}
q=\exp 2 \pi i \,\tau.
\end{equation}
The exponential of equation \eqref{eq:k3defnofmirrormap} is
\begin{equation}
q=\varphi \exp \left( \frac{h_1}{\varpi_0}\right),
\end{equation}
whose inversion gives us the $q$-expansion of $\varphi$. In fact, it is just the modular form $T(q)$ \eqref{eq:Beukerstwomodularforms} \cite{Beukers, BeukersPeters, Zagier}
\begin{equation} \label{eq:mirrormapmodularform}
\varphi=T(q).
\end{equation}
From formulas \eqref{eq:yukawaperiodsrepn} and \eqref{eq:k3defnofmirrormap}, the Yukawa coupling $\mathcal{Y}_{bp}$ can now be expressed as
\begin{equation} \label{eq:Yukawatauexpression}
\mathcal{Y}_{bp}=-2 C_0 \varpi_0^2\left( \frac{d\tau}{d\varphi} \right)^2.
\end{equation}
This formula immediately implies that the normalization of the Yukawa coupling $\mathcal{Y}_{bp} $, i.e.
\begin{equation} \label{eq:originalNormalizedYukawa}
\frac{1}{ \varpi_0^2 }\mathcal{Y}_{bp} \left( \frac{1}{2 \pi i}\frac{d \varphi}{d\tau} \right)^2,
\end{equation}
is just a nonzero constant. This property corresponds to the property that the Gromov-Witten invariants of K3 surfaces are trivial \cite{Nagura}.

\section{A new canonical period} \label{sec:varpi3}

In this section, we will construct a fourth order differential operator $\mathcal{D}$ from Ap\'ery's sequences, which also has the canonical periods $\varpi_i, i=0,1,2$ of $\mathcal{L}$ \eqref{eq:introPFoperatorK3} as its solutions. In addition, $\mathcal{D}$ has a fourth canonical solution $\varpi_3$, and the monodromy matrix of the canonical solutions $(\varpi_0,\varpi_1, \varpi_2,\varpi_3)$ is maximally unipotent at $\varphi=0$. The canonical solution $\varpi_3$ will play a crucial role in this paper.

First, let us construct a power series $\Pi_0(\phi)$ from the sequence $\{B_n \}_{n \geq 0}$ of Ap\'ery by
\begin{equation}
\Pi_0(\phi)=\sum_{n=0}^{\infty}B_n \phi^{n+1}=6\phi^2+\frac{351}{4}\phi^3+\frac{62531}{36} \phi^4+\frac{11424695}{288} \phi^5+\cdots.
\end{equation}
The recursion equation \eqref{eq:aperyrecursion} together with the initial condition $(B_0,B_1)=(0,6)$ immediately imply that $\Pi_0(\phi)$ is a solution of the fourth order differential operator
\begin{equation}\label{eq:varphiD}
\mathcal{D}= \left(1-34 \phi + \phi ^2\right) \theta^4-5 \left(1- 17\phi \right) \theta^3+3 (3-26\phi) \theta^2-(7-32\phi) \theta +(2-5\phi),~\theta =\phi \frac{d}{d\phi}.
\end{equation}
If we define a transformation between the variable $\varphi$ in Section \ref{sec:YukawaK3} and $\phi$ via,
\begin{equation}
\varphi=1/\phi,
\end{equation}
then there is a relation between the differential operator $\mathcal{D}$ \eqref{eq:varphiD} and the third order Picard-Fuchs operator $\mathcal{L}$ \eqref{eq:introPFoperatorK3}
\begin{equation}
\mathcal{D}=\varphi^{-2} \vartheta \cdot \mathcal{L},
\end{equation}
Hence the solutions $\{\varpi_i \}_{i=0}^2$ of the Picard-Fuchs operator $\mathcal{L}$ are also solutions of this new differential operator $\mathcal{D}$. 

Now let us find the fourth canonical solution of $\mathcal{D}$ by using the Frobenius method \cite{YangPeriods}. More explicitly, try a solution of the form
\begin{equation}
\mathcal{D}\,\left(\varphi^{\epsilon} \, \sum_{k=0}^{\infty} a_k(\epsilon)\, \varphi^k \right)=0,
\end{equation} 
where $\epsilon$ is a formal variable and $a_k(\epsilon)$ is a number. In order for this equation to be satisfied, to the lowest order, we must have
\begin{equation}
\epsilon^4=0,
\end{equation}
which is the indicial equation for $\mathcal{D}$. We have known that $\mathcal{D}$ already has three canonical solutions $\{\varpi_i \}_{i=0}^2$. Hence a fourth canonical solution $\varpi_3$ of $\mathcal{D}$, if exists, must be of the form
\begin{equation} \label{eq:varpi3formFrobenius}
\varpi_3(\varphi)=\frac{1}{(2 \pi i)^3}\left( \varpi_0(\varphi) \log^3 \varphi+3 h_1(\varphi) \log^2 \varphi+3 h_2(\varphi) \log \varphi+\sum_{n=1}^\infty c_{3,n} \varphi^n \right),
\end{equation}
where $h_1$ and $h_2$ are the two power series defined in the formulas \eqref{eq:h1defn} and \eqref{eq:h2defn}. Now plug this expression into $\mathcal{D} \varpi_3(\varphi)=0$, we obtain the recursion equation for $c_{3,n}$
\begin{equation}
\begin{aligned}
& (n+1)^4c_{3,n+1}-(34 n^4+85 n^3+78 n^2+32 n+5)c_{3,n} +n^3(n+1)c_{3,n-1} \\
&+12(n+1)^3c_{2,n+1}-3 \left(136 n^3+255 n^2+156 n+32\right) c_{2,n}+3 n^2 (4 n+3) c_{2,n-1} \\
&+36(n+1)^2c_{1,n+1} -18 \left(68 n^2+85 n+26\right) c_{1,n} +18 n (2 n+1) c_{1,n-1} \\
&+24(n+1)A_{n+1} -102(8n+5)A_n+6(4n+1)A_{n-1}=0
\end{aligned}
\end{equation}
with initial condition
\begin{equation}
c_{3,1}=-42,~c_{3,2}=-\frac{3033}{4}.
\end{equation}
For simplicity, let $h_3(\varphi)$ be
\begin{equation}
h_3(\varphi)=\sum_{n=1}^\infty c_{3,n} \varphi^n=-42 \varphi -\frac{3033}{4} \varphi^2 -\frac{522389}{36} \varphi^3 +\cdots.
\end{equation}

The solution $\varpi_3(\varphi)$ is indeed ``something new''! More precisely, we have
\begin{equation} \label{eq:varpi3inequality}
\varpi_0(\varphi) \varpi_3(\varphi) \neq \varpi_1(\varphi) \varpi_2(\varphi),
\end{equation}
thus the Picard-Fuchs operator $\mathcal{D}$ is not the symmetric cubic of the second order differential operator \eqref{eq:Lroot}. The four canonical solutions $\{\varpi_i \}_{i=0}^3$ form a basis for the solution space of $\mathcal{D}$. Define the column vector $\varpi$ by
\begin{equation}
\varpi= \left( \varpi_0,\varpi_1,\varpi_2,\varpi_3 \right)^\top.
\end{equation}
Then the monodromy matrix of $\varpi$ at $\varphi=0$ is
\begin{align}
\begin{pmatrix}
 1 & 0 & 0 & 0 \\
 1 & 1 & 0 & 0 \\
 1 & 2 & 1 & 0 \\
 1 & 3 & 3 & 1 \\
\end{pmatrix},
\end{align}
which is maximally unipotent \cite{MarkGross, KimYang}.

However, there are infinitely many fourth order differential operators which have the canonical solutions  $\varpi_i, i=0,1,2$ as its solutions, while also have a fourth canonical solution of the form \eqref{eq:varpi3formFrobenius}. For example, one such fourth order differential operator is 
\begin{equation} \label{eq:NewvarpiD}
\widetilde{\mathcal{D}}=\mathcal{D}+\varphi \mathcal{L},
\end{equation}
which also has a solution of the form \eqref{eq:varpi3formFrobenius}
\begin{equation} \label{eq:Newvarpi3formFrobenius}
\widetilde{\varpi}_3(\varphi)=\frac{1}{(2 \pi i)^3}\left( \varpi_0(\varphi) \log^3 \varphi+3 h_1(\varphi) \log^2 \varphi+3 h_2(\varphi) \log \varphi+ \tilde{h}_3 \right).
\end{equation}
Of course now $\tilde{h}_3$ is a different power series 
\begin{equation}
\tilde{h}_3=-48 \varphi -\frac{6765}{8} \varphi^2 -\frac{3507923}{216} \varphi ^3- \cdots.
\end{equation}
So the readers might be wondering why do we choose to study $\mathcal{D}$ in this paper? Does $\mathcal{D}$ have any interesting properties that make it the special one? This is exactly the subject of the next section.

\section{The Yukawa coupling and Instanton expansion} \label{sec:NewYukawa}

In this section, we will study the canonical solution $\varpi_3$ obtained in Section \ref{sec:varpi3}. We will construct a prepotential from $\varpi_3$ whose third derivative with respect to the mirror map $\tau$ \eqref{eq:k3defnofmirrormap} gives us a more interesting Yukawa coupling that is a weight-4 modular form. We will compute the instanton numbers for this Yukawa coupling, which are periodic integers.

The picture that we have obtained in the previous sections is very similar to that of the mirror symmetry of Calabi-Yau threefolds! We have four canonical solutions $\{\varpi_i \}_{i=0}^3$ for a fourth order differential operator $\mathcal{D}$ \eqref{eq:varphiD}, which are of the same type as that of the one-parameter families of Calabi-Yau threefolds, i.e. formula \eqref{eq:PeriodsCan}. Moreover, the quotient $\varpi_1/\varpi_0$ defines a mirror map $\tau$ \eqref{eq:k3defnofmirrormap}, which is similar to the mirror map $t$ \eqref{eq:mirrormap}. But of course, there is something different for $\mathcal{D}$!

From the definition of the mirror map $\tau$ in the formula \eqref{eq:k3defnofmirrormap}, we trivially have
\begin{equation} \label{eq:varpi1varpi0}
\varpi_1=\varpi_0 \tau.
\end{equation}
While the property that the Picard-Fuchs operator $\mathcal{L}$ \eqref{eq:introPFoperatorK3} is the symmetric square of the second order differential operator \eqref{eq:Lroot}, or equivalently formula \eqref{eq:varpi2defn}, tells us
\begin{equation} \label{eq:varpi2varpi0}
\varpi_2=\varpi_0 \tau^2.
\end{equation}
Notice that this is where the differences come in! Because of the inequality \eqref{eq:varpi3inequality}, the canonical solution $\varpi_3$ of $\mathcal{D}$ is not equal to $\varpi_0 \tau^3$. Instead, it is of the form
\begin{equation}
\varpi_3=\varpi_0 \tau^3 +\frac{1}{(2 \pi i)^3} \rho ,
\end{equation} 
where the power series $\rho$ is given by
\begin{equation}
\rho(\varphi) =h_3(\varphi)-\frac{h_1^3(\varphi)}{\varpi^2_0(\varphi)}.
\end{equation}
This power series can be considered as the correction obtained from the property that the operator $\mathcal{D}$ \eqref{eq:varphiD} is not the symmetric cubic of the second order differential operator \eqref{eq:Lroot}. 

Now, let us mimic the construction of the prepotential $\mathcal{F}$ \eqref{eq:prepotentialMirror} in the Calabi-Yau threefold case. Because of the formulas \eqref{eq:varpi1varpi0} and \eqref{eq:varpi2varpi0}, we have 
\begin{equation}
\frac{\varpi_1 \varpi_2}{\varpi_0^2} =\tau^3.
\end{equation}
So the only interesting term is the quotient $\varpi_3/\varpi_0$. Therefore, a natural definition of a prepotential for $\mathcal{D}$ \eqref{eq:varphiD} is  
\begin{equation} \label{eq:prePotentialK3Beukers}
\mathscr{F}_{\mathcal{D}}=\frac{\varpi_3}{\varpi_0}=\tau^3+\frac{1}{(2 \pi i)^3} \frac{ \rho (\varphi)}{\varpi_0(\varphi)}.
\end{equation}
While a natural definition of a Yukawa coupling $\mathscr{Y}_{\mathcal{D}}$ is
\begin{equation} \label{eq:YukawaNewK3}
\mathscr{Y}_{\mathcal{D}}=\frac{d^3\mathscr{F}_{\mathcal{D}}}{d\tau^3}.
\end{equation}
Because $q$ is $\exp (2 \pi i \tau)$, we have 
\begin{equation}
\frac{d}{d\tau} =(2 \pi i) q\frac{d}{dq}.
\end{equation}
Notice that there is no lose of generality in our definition which completely omits the terms involving $\varpi_1$ and $\varpi_2$, since they can only contribute to a cubic polynomial in $\tau$, whose third derivative is a constant.

Now plug the $q$-expansion of $\varphi$ \eqref{eq:mirrormapmodularform}, i.e. the modular form $T(q)$ \eqref{eq:Beukerstwomodularforms}, into $\mathscr{Y}_{\mathcal{D}}$, we obtain the $q$-expansion of $\mathscr{Y}_{\mathcal{D}}$ with the first several terms given by
\begin{equation}
\mathscr{Y}_{\mathcal{D}}=6 \left( 1-7q-59q^2-205 q^3-475 q^4-882 q^5- \cdots \right).
\end{equation}
Surprisingly, $\mathscr{Y}_{\mathcal{D}}$ is a weight-4 modular form that is the product of two weight-2 modular forms
\begin{equation}
\mathscr{Y}_{\mathcal{D}}=6 F(q)H(q).
\end{equation}
Here $F(q)$ is the weight-2 modular form in the formula \eqref{eq:Beukerstwomodularforms}. While $H(q)$ is the weight-2 modular form
\begin{equation}
H(q)=2 \Theta_{\text{hex}}^2(q^2)-\Theta_{\text{hex}}^2(q),
\end{equation}
where $\Theta_{\text{hex}}(q)$ is the theta series of the planar hexagonal lattice $A_2$ \cite{Conway, OEIStheta},
\begin{equation}
\Theta_{\text{hex}}(q)=\theta_3(0,\tau)\theta_3(0,3\tau)+\theta_2(0,\tau)\theta_2(0,3\tau).
\end{equation}

Following the definition of the instanton numbers in the mirror symmetry of Calabi-Yau threefolds, i.e. formula \eqref{eq:YukawaKahler}, we can define the instanton numbers for the Yukawa coupling $\mathscr{Y}_{\mathcal{D}}$ via the equation
\begin{equation} \label{eq:DefninstantNo}
\mathscr{Y}_{\mathcal{D}}=6 F(q)H(q)=6+\sum_{k=1}^\infty k^3 N_k \frac{q^k}{1-q^k}.
\end{equation}
We have computed the instanton numbers $N_k$ for small $k$, and the first six of them are given by
\begin{equation} \label{eq:K3InstantonNumber}
N_1 = -42, \, N_2=-39, \, N_3=-44, \, N_4=-39, \, N_5 =-42, \,N_6=-34.
\end{equation}
Very surprisingly, we have found that the instanton number $N_k$ is periodic
\begin{equation} \label{eq:BeukersK3Periodicity}
N_{k+6} = N_k.
\end{equation}
Since there is a factor $k^3$ in the definition of $N_k$ \eqref{eq:DefninstantNo}, so a priori we would expect $N_k$ is a fraction which has $k^3$ in its denominator. Thus the fact that $N_k$ is an integer is quite interesting. While the property that $N_k$ is periodic is certainly very surprising.

On the other hand, if we instead choose the operator $\widetilde{\mathcal{D}}$ \eqref{eq:NewvarpiD} from the beginning, and define the prepotential $\mathcal{F}_{\widetilde{\mathcal{D}}}$ by
\begin{equation} \label{eq:NewprePotentialK3Beukers}
\mathscr{F}_{\widetilde{\mathcal{D}}}=\frac{\widetilde{\varpi}_3}{\varpi_0},
\end{equation}
where the canonical solution $\widetilde{\varpi}_3$ is \eqref{eq:Newvarpi3formFrobenius}. Then we also have a new Yukawa coupling
\begin{equation} \label{eq:NewYukawaNewK3}
\mathscr{Y}_{\widetilde{\mathcal{D}}}=\frac{d^3\mathscr{F}_{\widetilde{\mathcal{D}}}}{d\tau^3}.
\end{equation}
However, there does not exist a nonzero integer $n$ such that the $q$-expansion of $n \mathscr{Y}_{\widetilde{\mathcal{D}}}$ has integral coefficients. So, we expect that the following properties pin down the choice of $\mathcal{D}$ uniquely
\begin{enumerate}
\item $\mathcal{D}$ is a fourth order differential operator which has $\{\varpi_i \}_{i=0}^2$ as its solutions;

\item There exists a fourth canonical solution for $\mathcal{D}$ that is of the form \eqref{eq:varpi3formFrobenius};

\item  The $q$-expansion of the Yukawa coupling $\mathscr{Y}_{\mathcal{D}}$ defined for $\mathcal{D}$  is a weight-4 modular form with integral coefficients (after a multiplication by a rational number if necessary).
\end{enumerate}
But we do not have a rigorous proof.

\section{The Dwork family of K3 surfaces } \label{sec:Dworkfamily}

In this section, we will generalize the results in Section \ref{sec:NewYukawa} to the Dwork family of K3 surfaces, which is the mirror family of the quartic K3 surfaces \cite{Hartmann, Nagura, YangK3}. We will obtain similar results to that of the previous section.

From the adjunction formula, a smooth quartic surface in $\mathbb{P}^3$ is a K3 surface. Its mirror family can be constructed from the Fermat pencil of K3 surfaces
\begin{equation} \label{eq:introFermatpencil}
\{ x_0^4+x_1^4+x_2^4+x_3^4-4\psi \,x_0x_1x_2x_3 =0 \} \subset \mathbb{P}^3
\end{equation} 
by first taking the quotient with respect to a $(\mathbb{Z}/4\mathbb{Z})^4$ action and then a minimal resolution of singularities \cite{Nagura}. Then we obtain a family of K3 surfaces over $\mathbb{P}^1$ with $\psi$ as its parameter, which is called the Dwork family of K3 surfaces \cite{Hartmann, Nagura, YangK3}. The nowhere-vanishing holomorphic twoform $\omega_\psi$ for the Dwork family satisfies a third order Picard-Fuchs equation $\mathcal{L}_{dk} \omega_\psi=0$, where the Picard-Fuchs operator $\mathcal{L}_{dk}$ is
\begin{equation} \label{eq:DworkPF}
\mathcal{L}_{dk}= \vartheta^3-4\varphi(4 \vartheta+1)(4 \vartheta+2)(4 \vartheta+3),~\text{with}~\varphi=1/(4\psi)^4~\text{and}~\vartheta=\varphi\,\frac{d}{d\varphi}.
\end{equation} 
It has two canonical solutions of the form \cite{Nagura}
\begin{equation} \label{eq:d3independentperiods}
\begin{aligned}
W_0&=\sum_{n=0}^\infty \frac{(4n)!}{(n!)^4}\varphi^n,\\
W_1&=\frac{1}{2 \pi i}\left(W_0 \cdot \log \varphi+4 \sum_{n=1}^\infty \frac{(4n)!}{(n!)^4}[\Psi(4n+1)-\Psi(n+1)] \varphi^n\right),\\
\end{aligned}
\end{equation}
where $\Psi(z)$ is the polygamma function
\begin{equation}
\Psi(z)=\frac{d}{dz} \log \Gamma(z).
\end{equation}
The Picard-Fuchs operator $\mathcal{L}_{dk} $ \eqref{eq:DworkPF} is also the symmetric square of a second order differential operator, hence there exists a third canonical solution $W_2$ for $\mathcal{L}_{dk}$ given by \cite{Nagura, YangK3}
\begin{equation}
W_2=W_1^2/W_0.
\end{equation}

The mirror map for the Dwork family of K3 surfaces is defined by the quotient
\begin{equation} \label{eq:quarticK3mirrormap}
\tau=\frac{W_1}{W_0}.
\end{equation}
Similarly, this equation can be inverted order by order, which gives us the $q$-expansion of $\varphi$ \cite{LianYau}
\begin{equation}\label{eq:quarticQexpansion}
\varphi=q-104q^2+6444q^3-311744q^4+ \cdots,~q=\exp (2 \pi i\,\tau).
\end{equation}
The three canonical solutions of $\mathcal{L}_{dk}$ can be expressed as
\begin{equation}
W_0,~W_0 \tau,~W_0\tau^2.
\end{equation}
Similarly, the fourth order differential operator 
\begin{equation}
\mathcal{D}_{dk}=\vartheta \mathcal{L}_{dk}
\end{equation}
has a fourth canonical solution $W_3$ of the form
\begin{equation}
W_3=W_0 \tau^3 +\frac{1}{(2 \pi i)^3}\rho_1(\varphi),
\end{equation}
where $\rho_1(\varphi)$ is a power series in $\varphi$. Define the prepotential $\mathscr{F}_{\mathcal{D}_{dk}}$ by
\begin{equation}
\mathscr{F}_{\mathcal{D}_{dk}}=\frac{W_3}{W_0}=\tau^3+\frac{1}{(2 \pi i)^3} \frac{\rho_1}{W_0}.
\end{equation}
The Yukawa coupling $\mathscr{Y}_{\mathcal{D}_{dk}}$ is similarly given by $d^3 \mathscr{F}_{\mathcal{D}_{dk}}/d\tau^3$, which has a $q$-expansion of the form
\begin{equation}
\mathscr{Y}_{\mathcal{D}_{dk}}=6-480 q-2400 q^2-13440q^3-17760q^4 - \cdots.
\end{equation}
But we do not know whether $\mathscr{Y}_{\mathcal{D}_{dk}}$ is a weight-4 modular form or not. Now let us define the instanton number $N_k$ for $\mathscr{Y}_{\mathcal{D}_{dk}}$ by the equation 
\begin{equation}
\mathscr{Y}_{\mathcal{D}_{dk}}=6+\sum_{k=1}^\infty k^3 N_k \frac{q^k}{1-q^k},
\end{equation}
from which we obtain
\begin{equation}
N_1=-480,~N_2=-240,~N_{k+2}=N_{k}.
\end{equation}
It is a very interesting question to find out whether these instanton numbers admit any geometric or physical interpretations.

\section{Conclusions and further prospects} \label{sec:furtherprospects}

In this paper, we have applied the ideas from the mirror symmetry of Calabi-Yau threefolds to study 
\begin{enumerate}
\item The two mysterious sequences constructed by Ap\'ery in his proof of the irrationality of $\zeta(3)$;

\item The modular forms and one-parameter family of K3 surfaces found by Beukers and Peters.
\end{enumerate}
The main result of this paper is that we have found a fourth order differential operator $\mathcal{D}$ \eqref{eq:varphiD}, and we have constructed a prepotential $\mathscr{F}_{\mathcal{D}}$ \eqref{eq:prePotentialK3Beukers}  from the four canonical solutions of $\mathcal{D}$ \eqref{eq:varphiD}. While the third derivative of $\mathscr{F}_{\mathcal{D}}$ with respect to the mirror map $\tau$ \eqref{eq:k3defnofmirrormap} gives us the Yukawa coupling $\mathscr{Y}_{\mathcal{D}}$ \eqref{eq:YukawaNewK3}. We have found that the $q$-expansion of the Yukawa coupling $\mathscr{Y}_{\mathcal{D}}$ is a weight-4 modular form. Moreover, we have defined instanton numbers for $\mathscr{Y}_{\mathcal{D}}$, which turn out to be integers. A very surprising result of this paper is that the instanton numbers for $\mathscr{Y}_{\mathcal{D}}$ are periodic, i.e. formula \eqref{eq:BeukersK3Periodicity}. We have also shown that our analysis can be generalized to the Dwork family of K3 surfaces, the mirror family for the quartic K3 surfaces, which gives us similar results.

There are many open questions left unanswered! Perhaps the most important open question is whether the computations in Section \ref{sec:NewYukawa} admit any interpretations in terms of the geometric information of K3 surfaces. For example, do the instanton numbers obtained in Section \ref{sec:NewYukawa} admit any interpretations in terms of certain invariants of K3 surfaces? What does the periodicity of the instanton numbers, i.e. formula \eqref{eq:BeukersK3Periodicity}, mean?

A more or less equivalent question is whether the prepotential $\mathscr{F}_{\mathcal{D}}$ \eqref{eq:prePotentialK3Beukers} and the Yukawa coupling $\mathscr{Y}_{\mathcal{D}}$ \eqref{eq:YukawaNewK3} have any relations with the string theory on K3 surfaces. Do they have any physical significance? Do the instanton numbers in Section \ref{sec:NewYukawa} admit any physical interpretations? Since the original normalized Yukawa coupling \eqref{eq:originalNormalizedYukawa} for the one-parameter family \eqref{eq:introK3family} is trivial, it would be very interesting if the Yukawa coupling \eqref{eq:YukawaNewK3} is a more meaningful K3-counterpart for the Yukawa coupling \eqref{eq:YukawaDefnKahler}.

Another related question is to properly understand why the differential operator $\mathcal{D}$ \eqref{eq:varphiD} is special? This has been intuitively discussed in the end of Section \ref{sec:NewYukawa}. But a more rigorous study on the properties of $\mathcal{D}$ will certainly help us understand the results in this paper better. It is also very interesting to see whether the results in this paper can be generalized to other one-parameter families of K3 surfaces. We have shown that the answer is yes for the Dwork family of K3 surfaces.



\begin{thebibliography}{99}






\bibitem{Apery}

R. Ap\'ery,  Irrationalit\'e de $\zeta(2)$ et $\zeta(3)$. Ast\'erisque. 61: 11–13.

\bibitem{Beukers}
F. Beukers, Irrationality proofs using modular forms. In Journ\'ees arithm\'etiques de Besan\c con (Besan\c con, 1985), Ast\'erisque 147–148 (1987), 271–283, 345.


\bibitem{BeukersPeters}
F. Beukers and C. Peters, A family of K3 surfaces and $\zeta(3)$. J. Reine Angew. Math. 351 (1984) 42–54.

\bibitem{PhilipXenia}

P. Candelas, X. C. de la Ossa, P. Green and L. Parkes, A Pair of Calabi-Yau Manifolds as an Exactly Soluble Superconformal Theory, Nuclear Physics B359 (1991) 21-74.


\bibitem{Conway}

J. H. Conway and N. J. A. Sloane, Sphere Packings, Lattices and Groups, Springer-Verlag, 1988, Section 6.2. Grundlehfen def mathematischen Wissenschaften 290.

\bibitem{CoxKatz}

D. Cox and S. Katz, Mirror Symmetry and Algebraic Geometry, American Mathematical Society. 

\bibitem{MarkGross}

M. Gross, D. Huybrechts and D. Joyce, Calabi-Yau Manifolds and Related Geometries. Springer.


\bibitem{Hartmann}

H. Hartmann, PERIOD- AND MIRROR-MAPS FOR THE QUARTIC K3. arXiv: 1101.4601.


\bibitem{KimYang}
M. Kim and W. Yang, Mirror symmetry, mixed motives and $\zeta(3)$. arXiv:1710.02344.

\bibitem{LianYau}

B. H. Lian and S. T. Yau, Arithmetic properties of mirror map and quantum coupling, Commun.
Math. Phys. 176 (1996) 163–91.

\bibitem{Nagura}

M. Nagura and K. Sugiyama, Mirror Symmetry of K3 and Torus. arXiv:9312159.

\bibitem{OEIS}

A186100, The On-Line Encyclopedia of Integer Sequences. \url{https://oeis.org/A186100}.

\bibitem{OEIStheta}

A004016, The On-Line Encyclopedia of Integer Sequences. \url{https://oeis.org/A004016}

\bibitem{YangK3}

K3 mirror symmetry, Legendre family and Deligne's conjecture for the Fermat quartic, Nuclear Physics B. \url{https://doi.org/10.1016/j.nuclphysb.2020.115303}.  arXiv: 2004.00820.

\bibitem{YangPeriods}

W. Yang, Periods of CY $n$-folds and mixed Tate motives, a numerical study. arXiv:1908.09965

\bibitem{Yang1}

W. Yang, Double zeta values and Picard-Fuchs equation. arXiv:1910.09576


\bibitem{Zagier}

D. Zagier, Arithmetic and topology of differential equations, in Proceedings of the Seventh European
Congress of Mathematics (Berlin, July 18–22, 2016), Editors V. Mehrmann, M. Skutella, European Mathematical Society, Berlin, 2018, 717–776.







\end{thebibliography}
\end{document}